\documentclass[preprint,11pt,times]{elsarticle}
\usepackage[T1]{fontenc}
\usepackage{amssymb,amsmath,amsthm,mathtools}
\usepackage{geometry}
\usepackage{setspace}
\usepackage{enumitem,tikz-cd}
\usepackage{needspace}

\newtheorem{theorem}{Theorem}[section]
\newtheorem{corollary}[theorem]{Corollary}
\newtheorem{lemma}[theorem]{Lemma}
\newtheorem{proposition}[theorem]{Proposition}
\theoremstyle{definition}

\newtheorem{criterion}[theorem]{Criterion}
\theoremstyle{remark}

\newtheorem{example}[theorem]{Example}
\numberwithin{equation}{section}
\usepackage[colorlinks,linkcolor=blue,anchorcolor=blue,citecolor=blue]{hyperref}
\usepackage{color}
\hypersetup{pdftitle={Linear Dynamics of Composition Operators on Zygmund Spaces of Rooted Trees},pdfauthor={Xiang Chen, Siwei Song, Liang Zhang, Zehua Zhou}}

\newcommand{\N}{\mathbb N}
\newcommand{\Nzero}{\mathbb N_0}
\newcommand{\C}{\mathbb C}
\newcommand{\F}{\mathcal F}
\newcommand{\Z}{\mathcal Z}
\newcommand{\Lip}{\mathcal L}
\newcommand{\ind}{\mathbf1}
\newcommand{\supp}{\operatorname{supp}}
\newcommand{\Per}{\operatorname{Per}}

\newcommand{\norm}[1]{\lVert#1\rVert}

\begin{document}
\begin{frontmatter}
\title{Linear Dynamics of Composition Operators\\
on Zygmund Spaces of Rooted Trees}
\author[aff1]{Xiang Chen\corref{cor1}}
\ead{2020233018@tju.edu.cn}

\author[aff1]{Siwei Song}
\ead{2025233007@tju.edu.cn}

\author[aff2]{Liang Zhang}
\ead{168zhangliang2011@163.com}

\author[aff1]{Zehua Zhou}
\ead{zehuazhoumath@aliyun.com; zhzhou@tju.edu.cn}

\cortext[cor1]{Corresponding author}

\affiliation[aff1]{organization={School of Mathematics, Tianjin University},
  city={Tianjin}, postcode={300350},
  country={P.R. China}}

\affiliation[aff2]{organization={School of Marine Science and Technology, Tianjin University},
  city={Tianjin}, postcode={300050},
  country={P.R. China}}
\begin{abstract}
We conjugate composition operators on the Zygmund and little Zygmund spaces of an infinite, locally finite rooted tree~$T$ by the second-difference isometries. A path Green kernel gives row-finite matrices that characterize boundedness and yield exact operator norms. A rigidity result holds on the chain~$\Nzero$. Boundedness on $\Z_0(\Nzero)$ is equivalent to $\varphi(n)=qn+b$ for all sufficiently large~$n$, with $q\in\Nzero$ and $b\in\mathbb Z$. For these symbols, $\varphi(n)>n$ eventually characterizes Li--Yorke chaos of $C_\varphi$ on either Zygmund space, allowing noninjectivity and periodic vertices. If $C_\varphi$ is bounded on $\Z_0(T)$ and $\varphi^k(T)$ eventually avoids every finite set, then for each $\lambda\ne0$, $\lambda C_\varphi$ is Li--Yorke chaotic on either space exactly when the norms of its powers are unbounded. For $\varphi_q(n)=q(n+1)-1$ with integer~$q\ge2$, the corresponding conditions are $|\lambda|q^2>1$ on both Zygmund spaces and $|\lambda|q>1$ on the little Lipschitz space. For bounded $C_\varphi$ on $\Z_0(\Nzero)$, topological mixing and Devaney chaos each hold exactly for injective symbols with no periodic vertices and eventual slope $q\ge2$. For such bounded operators with strictly increasing $\varphi$ and $\varphi(0)>0$, $q\ge2$ also characterizes hypercyclicity and frequent hypercyclicity.
\end{abstract}
\end{frontmatter}

\section{Introduction}
A self-map $\varphi$ of an infinite, locally finite rooted tree $T$
induces the composition operator $C_\varphi f=f\circ\varphi$ on
$\F(T)=\C^T$. We investigate how the passage from first to second
differences affects the boundedness and linear dynamics of this
operator.

Colonna and Easley~\cite{CE} studied multiplication on the Lipschitz
space of a tree, and Allen, Colonna, and Easley~\cite{ACE}
characterized bounded composition operators on that space.
Colonna and Mart\'inez-Avenda\~no~\cite{CM} examined boundedness,
spectra, and hypercyclicity on the little Lipschitz space.
L\'opez-Mart\'inez~\cite{LM} used first-difference isometries to
compute exact norms. On the little Lipschitz space of the chain
$\Nzero$, he showed that a bounded $C_\varphi$ is hypercyclic if and
only if it satisfies the Frequent Hypercyclicity Criterion, provided
$\varphi$ is strictly increasing and $\varphi(0)>0$.
Related studies treat discrete differentiation and
cancellation~\cite{AJ,LiuKong}, shifts on tree Lipschitz and Hardy
spaces~\cite{AMHardy,Kawamura,MR,ZZ}, weighted shifts on directed
trees~\cite{AA,GPdyn,MA}, and composition operators on discrete
weighted spaces~\cite{ACMP,ZCCZ}.

Delgadillo-Garc\'ia and Mart\'inez-Avenda\~no~\cite{DGMA} studied
forward and backward shifts on the Zygmund spaces of a tree.
We adopt their discrete difference $D$. The spaces $\Z(T)$ and
$\Z_0(T)$ consist of functions for which $D^2f$ is bounded or tends
to zero with depth, respectively, and both have norm
$\norm{f}_{\Z}=\norm{D^2f}_\infty$. An edge increment of $f\circ\varphi$
is a signed sum of first differences of $f$ along the path joining
the images of the endpoints. At depth at least two, the second difference
involves the images of a vertex, its parent, and its grandparent.
The corresponding weighted contributions along root paths can cancel.
For instance, $\varphi(n)=\lfloor n/2\rfloor$ induces a bounded
operator on $\Lip_0(\Nzero)$ and takes $f(n)=n+1$ outside
$\Z_0(\Nzero)$ (Example~\ref{ex:floor}).

\Needspace{6\baselineskip}
The second difference identifies $\Z(T)$ and $\Z_0(T)$ isometrically
with $\ell^\infty(T)$ and $c_0(T)$, respectively. Its inverse is
$P^2$, where $P$ sums along root paths. With $|a|$ the depth of $a$
and $w\preceq a$ denoting membership in the root path to $a$, the
reconstruction formula of~\cite[Lemma 5.2]{DGMA} takes the form
\[
G(a,w)=(|a|-|w|+1)\ind_{\{w\preceq a\}},
\qquad
(P^2x)(a)=\sum_{w\in T}G(a,w)x(w).
\]
Since $D^2G(\cdot,w)=e_w$, the point mass at $w$, we call $G$ the
\emph{path Green kernel}. In the row-finite matrix of $D^2C_\varphi P^2$, each entry combines
the path contributions before the absolute row sums are formed. A uniform bound on the absolute row sums characterizes boundedness
on $\Z(T)$. On $\Z_0(T)$, boundedness requires the same row bound and
convergence of each column to zero. Since the coefficients are integers,
the latter condition is equivalent to finite support of every column.
In either case, the supremum of the row sums gives the operator norm.
On the chain, $G(n,0)=n+1$, so finite support of the root column
forces the symbol's second differences to vanish eventually.
Conversely, an affine tail satisfies both matrix conditions
(Theorem~\ref{thm:affine}).

To study Li--Yorke chaos, we use $C_\varphi^k=C_{\varphi^k}$ to compute
power norms, as in the norm approach to weighted shifts~\cite{KJCK}.
When $\varphi^k(T)$ eventually avoids every finite subset of $T$,
$C_\varphi^kf$ is eventually zero for every function $f$ with finite
support on $T$. For bounded $C_\varphi$ on
$\Z_0(T)$, these functions form a dense stable set there. This fact
and the common power norms yield the Li--Yorke criterion for nonzero
scalar multiples on both spaces (Theorem~\ref{thm:LY-tree}). On
$\Nzero$, with boundedness on the little space, eventual forward
displacement $\varphi(n)>n$ characterizes Li--Yorke chaos of
$C_\varphi$ itself. For sufficiency we work with functions vanishing
at periodic vertices (Theorem~\ref{thm:LY-chain}). For $\varphi(n)=n+1$, the composition operator exhibits Li--Yorke
chaos and has no dense orbit on the little space. The expanding affine
family gives distinct scalar thresholds at first and second order (Theorem~\ref{thm:LY-threshold}).

Since $\Z(T)$ is nonseparable and $\Z_0(T)$ is separable, we consider
dense-orbit properties in the little space. On the chain, for bounded
$C_\varphi$, topological mixing and Devaney chaos each hold exactly when
the symbol is injective, has no periodic vertices, and has eventual
slope $q\ge2$ (Theorem~\ref{thm:chaos}). The sufficiency proofs use
an interpolating right inverse with summable iterates on finitely
supported functions. Strict increase and $\varphi(0)>0$ give an injective symbol with
$\varphi(n)>n$ for every $n$. Under these additional hypotheses the slope bound also characterizes
hypercyclicity, frequent hypercyclicity, and satisfaction of the
Hypercyclicity Criterion~\cite[Chapter 3]{LinearChaos} (Theorem~\ref{thm:frequent} and
Corollary~\ref{cor:increasing-chaos}).

Section~\ref{sec:conjugacy} develops the isometries, the path Green
kernel, and the finite-support approximation.
Section~\ref{sec:boundedness} gives the boundedness criteria,
exact norm formulas, and the chain characterization.
Section~\ref{sec:dynamics} treats Li--Yorke chaos and the
dense-orbit properties.

\section{Preliminaries and kernel representation}\label{sec:conjugacy}
\subsection{Function spaces and conjugacy}\label{subsec:function-spaces}
Let $T$ be an infinite, locally finite rooted tree with root $o$. We identify the tree with its vertex set and allow leaves. Each nonroot vertex $v$ has a unique parent $v^-$. For $|v|\ge2$, put $v^{--}=(v^-)^-$. The depth is $|v|=d(o,v)$, where $d$ is the edge-counting distance. We write $w\preceq v$ when $w$ lies on the root path $[o,v]$, including its endpoints, and set
\[
S_w=\{v\in T:w\preceq v\}.
\]
Local finiteness implies that every depth ball is finite; in particular, $T$ is countably infinite. Let $\F(T)=\C^T$.

The Lipschitz-space representation of~\cite[Remarks 2.1 and 2.2]{LM} uses the algebraic difference operator $D:\F(T)\to\F(T)$ given by
\begin{equation}\label{eq:D}
(Df)(o)=f(o),\qquad (Df)(v)=f(v)-f(v^-)\quad(v\ne o).
\end{equation}
The sequence spaces
\begin{align*}
\ell^\infty(T)&=\{x\in\F(T):\sup_{v\in T}|x(v)|<\infty\},\\
c_0(T)&=\{x\in\F(T):x(v)\to0\text{ as }|v|\to\infty\},
\end{align*}
carry the supremum norm. The Lipschitz and little Lipschitz spaces are
\begin{equation}\label{eq:Lspaces}
\Lip(T)=\{f\in\F(T):Df\in\ell^\infty(T)\},\qquad
\Lip_0(T)=\{f\in\F(T):Df\in c_0(T)\},
\end{equation}
with
\[
\norm f_{\Lip}=\norm{Df}_\infty
=\max\bigl\{|f(o)|,\sup_{v\ne o}|f(v)-f(v^-)|\bigr\}.
\]
We use the maximum norm of~\cite{CM,LM}, which is equivalent to the sum norm in~\cite{ACE}. All exact norm formulas below refer to this convention.

Define the path-sum operator $P:\F(T)\to\F(T)$ by
\begin{equation}\label{eq:P}
(Px)(a)=\sum_{w\preceq a}x(w).
\end{equation}
Adjacent path sums give $DP=I$, including at the root, while telescoping gives $PD=I$. Hence the restrictions
\begin{equation}\label{eq:first-restrictions}
D_{\Lip}:=D|_{\Lip(T)}:\Lip(T)\to\ell^\infty(T),\qquad
D_{\Lip,0}:=D|_{\Lip_0(T)}:\Lip_0(T)\to c_0(T)
\end{equation}
are surjective linear isometries, with inverses $P|_{\ell^\infty(T)}$ and $P|_{c_0(T)}$, respectively.

For $a,b\in T$, telescoping gives the standard point-value and path estimates; see~\cite[Lemma 2.1]{ACE}:
\begin{equation}\label{eq:L-estimates}
|f(a)|\le(|a|+1)\norm f_{\Lip},\qquad
|f(a)-f(b)|\le d(a,b)\norm f_{\Lip}
\quad(f\in\Lip(T)).
\end{equation}

For $e_w(v)=\delta_{w,v}$, one has~\cite[Remark 2.1]{LM}
\[
Pe_w=\ind_{S_w},\qquad D\ind_{S_w}=e_w.
\]
Thus the sector indicators span a dense subspace of $\Lip_0(T)$. Finitely supported functions are also dense~\cite{CM}: linearly cutting off a sector indicator over $m\ge1$ further levels gives finite support and $\Lip$-error at most $1/m$. This argument uses only local finiteness. The second-order analogue on the full tree is proved in Lemma~\ref{lem:finite-density}.

For the composition operator $C_\varphi$ on $\F(T)$, put $B_\varphi=DC_\varphi P$. If $C_\varphi\in\mathcal B(\Lip_0(T))$, then
\[
B_\varphi|_{c_0(T)}
=D_{\Lip,0}\,(C_\varphi|_{\Lip_0(T)})\,D_{\Lip,0}^{-1},
\]
as expressed by the commutative diagram
\begin{equation}\label{eq:L-diagram}
\begin{tikzcd}[column sep=4em,row sep=2.2em]
\Lip_0(T)\arrow[r,"C_\varphi|_{\Lip_0(T)}"]\arrow[d,"D_{\Lip,0}"']
&\Lip_0(T)\arrow[d,"D_{\Lip,0}"]\\
c_0(T)\arrow[r,"B_\varphi|_{c_0(T)}"']&c_0(T).
\end{tikzcd}
\end{equation}
The analogous conjugacy holds for $\Lip(T)$ and $\ell^\infty(T)$ with $D_{\Lip}$ under the corresponding boundedness hypothesis. These conjugacies preserve the dynamical properties considered below~\cite[Remark 2.2]{LM}.

For comparison, the classical criterion~\cite[Theorem 3.2]{ACE} and exact norm formula~\cite[Corollary 2.3]{LM} state that $C_\varphi$ is bounded on $\Lip(T)$ exactly when the following quantity is finite:
\begin{equation}\label{eq:L-norm}
\norm{C_\varphi}_{\Lip\to\Lip}
=\max\left\{1+|\varphi(o)|,\ \sup_{v\ne o}d\bigl(\varphi(v),\varphi(v^-)\bigr)\right\}.
\end{equation}
The norm is the same on $\Lip_0(T)$ when that space is preserved; the finite bound alone need not ensure this preservation~\cite{CM,LM}.

Applying $D$ twice on $\F(T)$ gives the second difference $D^2=D\circ D$. With the boundary convention in~\cite[Definition 2.1]{DGMA}, it has the form
\begin{equation}\label{eq:D2}
(D^2f)(v)=
\begin{cases}
f(o),&v=o,\\
f(v)-2f(o),&|v|=1,\\
f(v)-2f(v^-)+f(v^{--}),&|v|\ge2.
\end{cases}
\end{equation}
The Zygmund and little Zygmund spaces are
\begin{equation}\label{eq:spaces}
\Z(T)=\{f\in\F(T):D^2f\in\ell^\infty(T)\},\qquad
\Z_0(T)=\{f\in\F(T):D^2f\in c_0(T)\}.
\end{equation}
Both spaces are endowed with
\begin{equation}\label{eq:znorm}
\norm f_{\Z}=\norm{D^2f}_\infty.
\end{equation}
The continuous inclusions $\Lip(T)\subset\Z(T)$ and $\Lip_0(T)\subset\Z_0(T)$ follow from $\norm f_{\Z}\le2\norm f_{\Lip}$, since $D$ is bounded on $\ell^\infty(T)$ and preserves $c_0(T)$; compare~\cite[Section 2]{DGMA}.
By definition, $D$ sends $\Z(T)$ into $\Lip(T)$ and $\Z_0(T)$ into $\Lip_0(T)$. Write
\begin{equation}\label{eq:second-restrictions}
D_{\Z}:=D|_{\Z(T)}:\Z(T)\to\Lip(T),\qquad
D_{\Z,0}:=D|_{\Z_0(T)}:\Z_0(T)\to\Lip_0(T).
\end{equation}

\begin{proposition}\label{prop:isometry}
The maps $D_{\Z}$ and $D_{\Z,0}$ are surjective linear isometries, with inverses $P|_{\Lip(T)}$ and $P|_{\Lip_0(T)}$. Consequently,
\begin{align*}
J&:=D_{\Lip}\circ D_{\Z}:\Z(T)\to\ell^\infty(T),& Jf&=D^2f,\\
J_0&:=D_{\Lip,0}\circ D_{\Z,0}:\Z_0(T)\to c_0(T),& J_0f&=D^2f
\end{align*}
are surjective linear isometries. Their inverses are the corresponding restrictions of $P^2$, given by
\begin{equation}\label{eq:inverse-path}
(P^2x)(a)=\sum_{w\preceq a}(|a|-|w|+1)x(w).
\end{equation}
\end{proposition}
\begin{proof}
The identity $D^2(Pg)=Dg$ shows that $P$ sends $\Lip(T)$ into $\Z(T)$ and $\Lip_0(T)$ into $\Z_0(T)$. Thus $DP=PD=I$ and $\norm{Df}_{\Lip}=\norm f_{\Z}$ give the two inverse isometries. Composing with~\eqref{eq:first-restrictions} gives $J$ and $J_0$. In
\[
(P^2x)(a)=\sum_{u\preceq a}\sum_{w\preceq u}x(w),
\]
each $x(w)$ with $w\preceq a$ occurs $|a|-|w|+1$ times, once for each vertex between $w$ and $a$. This yields~\eqref{eq:inverse-path}.
\end{proof}

The two applications of $D$ are thus represented by the commutative diagrams
\begin{equation}\label{eq:D-factorization}
\begin{tikzcd}[column sep=2.8em]
\Z(T)\arrow[r,"D_{\Z}"']\arrow[rr,bend left=28,"J"]
&\Lip(T)\arrow[r,"D_{\Lip}"']&\ell^\infty(T)
\end{tikzcd}
\qquad
\begin{tikzcd}[column sep=2.8em]
\Z_0(T)\arrow[r,"D_{\Z,0}"']\arrow[rr,bend left=28,"J_0"]
&\Lip_0(T)\arrow[r,"D_{\Lip,0}"']&c_0(T).
\end{tikzcd}
\end{equation}
The four maps in~\eqref{eq:D-factorization} are restrictions of the same algebraic operator~\eqref{eq:D}. On $\F(T)$ one has the self-composition $D^2=D\circ D$; on the normed spaces the corresponding compositions are $J=D_{\Lip}\circ D_{\Z}$ and $J_0=D_{\Lip,0}\circ D_{\Z,0}$. The inverse composition $P^2$ gives the path reconstruction formula $f=P^2D^2f$ from~\cite[Lemma 5.2]{DGMA}.

\begin{corollary}\label{cor:structure}
The space $\Z(T)$ is a nonseparable Banach space, and $\Z_0(T)$ is a separable closed subspace of $\Z(T)$. Moreover,
\begin{equation}\label{eq:pointbounds}
|Df(a)|\le(|a|+1)\norm f_{\Z},\qquad
|f(a)|\le\frac{(|a|+1)(|a|+2)}2\norm f_{\Z}.
\end{equation}
In particular, all point evaluations are continuous on both spaces.
\end{corollary}
\begin{proof}
The assertions about the spaces follow from Proposition~\ref{prop:isometry} and the corresponding properties of $\ell^\infty(T)$ and $c_0(T)$. Taking absolute values in the path sums for $Df=P(D^2f)$ and $f=P^2(D^2f)$ gives~\eqref{eq:pointbounds}.
\end{proof}

For the composition operator defined above, introduce the algebraic conjugate
\begin{equation}\label{eq:algebraicA}
A_\varphi=D^2C_\varphi P^2:\F(T)\to\F(T).
\end{equation}
Since $D^2$ is invertible on $\F(T)$,
\begin{equation}\label{eq:intertwine}
D^2C_\varphi=A_\varphi D^2,\qquad
C_\varphi=P^2A_\varphi D^2.
\end{equation}

\begin{proposition}\label{prop:conjugacy}
For $X=\Z(T)$ or $X=\Z_0(T)$, preservation $C_\varphi(X)\subseteq X$ is equivalent to boundedness on $X$. Furthermore,
\begin{align}
C_\varphi\in\mathcal B(\Z(T))
&\iff A_\varphi|_{\ell^\infty(T)}\in\mathcal B(\ell^\infty(T)),\label{eq:big-conjugacy}\\
C_\varphi\in\mathcal B(\Z_0(T))
&\iff A_\varphi|_{c_0(T)}\in\mathcal B(c_0(T)).\label{eq:small-conjugacy}
\end{align}
Under the respective boundedness hypotheses,
\[
A_\varphi|_{\ell^\infty(T)}=J(C_\varphi|_{\Z(T)})J^{-1},\qquad
A_\varphi|_{c_0(T)}=J_0(C_\varphi|_{\Z_0(T)})J_0^{-1},
\]
and each conjugate has the same operator norm as $C_\varphi$ on the corresponding Zygmund space.
\end{proposition}
\begin{proof}
Assume $C_\varphi(X)\subseteq X$ and let $f_n\to f$ and $C_\varphi f_n\to g$ in $X$. Continuity of point evaluations gives $g(v)=f(\varphi(v))$ for each $v$, so the closed graph theorem gives boundedness. Restricting~\eqref{eq:intertwine} to the two spaces via the isometries $J$ and $J_0$ gives the conjugacy and norm identities.
\end{proof}

The corresponding commutative diagrams are
\begin{equation}\label{eq:Z-diagrams}
\begin{tikzcd}[column sep=3em,row sep=2.2em]
\Z(T)\arrow[r,"C_\varphi|_{\Z(T)}"]\arrow[d,"J"']
&\Z(T)\arrow[d,"J"]\\
\ell^\infty(T)\arrow[r,"A_\varphi|_{\ell^\infty(T)}"']&\ell^\infty(T)
\end{tikzcd}
\qquad
\begin{tikzcd}[column sep=3em,row sep=2.2em]
\Z_0(T)\arrow[r,"C_\varphi|_{\Z_0(T)}"]\arrow[d,"J_0"']
&\Z_0(T)\arrow[d,"J_0"]\\
c_0(T)\arrow[r,"A_\varphi|_{c_0(T)}"']&c_0(T).
\end{tikzcd}
\end{equation}

\subsection{The path Green kernel}\label{subsec:kernel}
The coefficients of $P^2$ are supported on finite root paths. For $a,w\in T$, put
\begin{equation}\label{eq:HG}
H(a,w)=\ind_{\{w\preceq a\}},\qquad
G(a,w)=(|a|-|w|+1)\ind_{\{w\preceq a\}}.
\end{equation}
For a fixed $w$, put $H_w(a)=H(a,w)$ and $G_w(a)=G(a,w)$. Thus $H_w=Pe_w=\ind_{S_w}$ is the sector function from the first-order representation. We call $G$ the \emph{path Green kernel} of $D^2$. Its weights are those of the reconstruction formula~\eqref{eq:inverse-path}.

\Needspace{24\baselineskip}
\begin{lemma}\label{lem:kernel}
For all $a,w\in T$,
\begin{enumerate}[label=\textup{(\roman*)}]
\item $DG_w=H_w$, $DH_w=e_w$, and $D^2G_w=e_w$;
\item $G_w$ is the unique solution of $D^2u=e_w$ in $\F(T)$;
\item for every $x,f\in\F(T)$,
\begin{equation}\label{eq:G-reconstruction}
(P^2x)(a)=\sum_{w\in T}G(a,w)x(w),\qquad
f(a)=\sum_{w\in T}G(a,w)(D^2f)(w);
\end{equation}
\item $G\ge0$, $\supp_wG(a,w)=[o,a]$, $\supp_aG(a,w)=S_w$, and
\begin{equation}\label{eq:G-row}
G(w,w)=1,\qquad \sum_{w\in T}G(a,w)=\frac{(|a|+1)(|a|+2)}2;
\end{equation}
\item the kernel-product formula is
\begin{equation}\label{eq:G-product}
G(a,w)=\sum_{u\in T}H(a,u)H(u,w).
\end{equation}
All sums are finite for fixed output vertex $a$.
\end{enumerate}
\end{lemma}
\begin{proof}
Part (iii) follows from~\eqref{eq:inverse-path} and $P^2D^2=I$. Taking $x=e_w$ gives $G_w=P^2e_w$. Since $H_w=Pe_w$ and $D=P^{-1}$, (i) and (ii) follow. Positivity, the support statements, and $G(w,w)=1$ follow from~\eqref{eq:HG}. Summing the weights on $[o,a]$ gives~\eqref{eq:G-row}.

For (v), $H(a,u)H(u,w)=1$ exactly when $w\preceq u\preceq a$. Thus
\[
\sum_{u\in T}H(a,u)H(u,w)
=\#\{u\in T:w\preceq u\preceq a\}
=G(a,w),
\]
by the same path count as in the proof of Proposition~\ref{prop:isometry}.
\end{proof}

\begin{corollary}\label{cor:eval-norm}
For $\varepsilon_a(f)=f(a)$,
\[
\norm{\varepsilon_a}_{\Z(T)^*}=\norm{\varepsilon_a}_{\Z_0(T)^*}
=\frac{(|a|+1)(|a|+2)}2.
\]
\end{corollary}
\begin{proof}
The bound in~\eqref{eq:pointbounds} is attained at $P^2\ind_{[o,a]}\in\Z_0(T)$. Its second difference has supremum norm one, and its value at $a$ is the row sum in~\eqref{eq:G-row}.
\end{proof}

Applying the second difference to the kernel representation gives the coefficients of $A_\varphi$.

\begin{proposition}\label{prop:matrix}
The operator in~\eqref{eq:algebraicA} has the row-finite representation
\begin{equation}\label{eq:A-matrix}
(A_\varphi x)(v)=\sum_{w\in T}a_\varphi(v,w)x(w),
\end{equation}
where
\begin{equation}\label{eq:A-coeff}
a_\varphi(v,w)=
\begin{cases}
G(\varphi(o),w),&v=o,\\
G(\varphi(v),w)-2G(\varphi(o),w),&|v|=1,\\
G(\varphi(v),w)-2G(\varphi(v^-),w)+G(\varphi(v^{--}),w),&|v|\ge2.
\end{cases}
\end{equation}
Every coefficient is a real integer. For each $w\in T$,
\begin{equation}\label{eq:column-test}
G_w=J_0^{-1}e_w\in\Z_0(T),\quad \norm{G_w}_{\Z}=1,
\qquad (D^2C_\varphi G_w)(v)=a_\varphi(v,w).
\end{equation}
The span of $\{G_w:w\in T\}$ is dense in $\Z_0(T)$.
\end{proposition}
\begin{proof}
Substituting~\eqref{eq:G-reconstruction} into $D^2C_\varphi P^2$ and using~\eqref{eq:D2} gives~\eqref{eq:A-matrix} and~\eqref{eq:A-coeff}. Each row is an integer combination of at most three integer-valued Green rows, each supported on a finite root path. Since $J_0G_w=e_w$ and $A_\varphi e_w=D^2C_\varphi G_w$, \eqref{eq:column-test} follows. The span of the $G_w$ is dense because $J_0$ maps it onto the finitely supported sequences.
\end{proof}

For the dynamics, finite support is needed in the original function values, not only in the second differences. The following approximation also preserves prescribed zero values.

\begin{lemma}\label{lem:finite-density}
Let $\F_{00}(T)=c_{00}(T)$ be the finitely supported functions on $T$. Then $\F_{00}(T)$ is dense in $\Z_0(T)$. For every finite set $F\subseteq T$, put
\[
Y_F=\{f\in\Z_0(T):f(a)=0\text{ for all }a\in F\}.
\]
The subspace $\F_{00}(T)\cap Y_F$ is dense in $Y_F$.
\end{lemma}
\begin{proof}
By~\eqref{eq:D2}, $D^2$ can enlarge support by at most two descendant levels. Local finiteness therefore gives $\F_{00}(T)\subseteq\Z_0(T)$. By Proposition~\ref{prop:matrix}, it suffices to approximate each $G_w$. Choose $\theta\in C^2([0,\infty))$ equal to $1$ on $[0,1]$ and to $0$ on $[2,\infty)$. For an integer $R\ge2$, set
\[
h_{w,R}(v)=
\begin{cases}
(t+1)\theta(t/R),&w\preceq v,\quad t=|v|-|w|,\\
0,&w\not\preceq v.
\end{cases}
\]
Its support lies in the finite depth ball $\{v:|v|\le|w|+2R\}$. At the root and first level, $h_{w,R}$ agrees with $G_w$, so their boundary differences in~\eqref{eq:D2} agree. Outside $S_w$, every ancestor contributing to the second difference is also outside $S_w$. Within $S_w$, the functions agree at every vertex of relative depth $t\le R$ and at its ancestors. Thus the second-difference error vanishes there.

For $t>R$, the second difference of $G_w$ is zero. The values of $h_{w,R}$ along the parent path are those of $s\mapsto(s+1)\theta(s/R)$ at $t,t-1,t-2$, regardless of branching. Twice applying the fundamental theorem of calculus bounds this second difference by the supremum of the absolute value of
\[
\frac{2}{R}\theta'(s/R)+\frac{s+1}{R^2}\theta''(s/R).
\]
The derivatives vanish outside $R\le s\le2R$. On this interval, $(s+1)/R^2\le3/R$, so the bound is uniform over all vertices in the transition region. For $t>2R+2$, all three cutoff values vanish. Consequently,
\begin{equation}\label{eq:Green-cutoff}
\norm{h_{w,R}-G_w}_{\Z}
\le\frac{2\norm{\theta'}_\infty+3\norm{\theta''}_\infty}{R}
\longrightarrow0.
\end{equation}
For $f\in Y_F$, take finitely supported $g_j\to f$ in $\Z_0(T)$. Then $g_j-\sum_{a\in F}g_j(a)e_a$ belongs to $\F_{00}(T)\cap Y_F$ and converges to $f$ because $F$ is finite and point evaluations are continuous.
\end{proof}

\section{Boundedness and exact norms}\label{sec:boundedness}
The absolute row sums of the conjugate matrix determine its norm, and the columns determine whether a bounded matrix preserves $c_0(T)$.

For the coefficients in~\eqref{eq:A-coeff}, put
\begin{equation}\label{eq:M}
M_\varphi=\sup_{v\in T}\sum_{w\in T}|a_\varphi(v,w)|.
\end{equation}
Each row sum is a finite nonnegative integer.

\Needspace{9\baselineskip}
\begin{theorem}\label{thm:boundedness}
Let $\varphi:T\to T$ be arbitrary.
\begin{enumerate}[label=\textup{(\alph*)}]
\item The operator $C_\varphi$ is bounded on $\Z(T)$ if and only if $M_\varphi<\infty$. In this case
\begin{equation}\label{eq:big-norm}
\norm{C_\varphi}_{\Z\to\Z}=M_\varphi.
\end{equation}
\item The operator $C_\varphi$ is bounded on $\Z_0(T)$ if and only if $M_\varphi<\infty$ and, for every $w\in T$,
\begin{equation}\label{eq:finite-columns}
\{v\in T:a_\varphi(v,w)\ne0\}\text{ is finite}.
\end{equation}
Whenever these conditions hold, it is also bounded on $\Z(T)$ and
\begin{equation}\label{eq:small-norm}
\norm{C_\varphi}_{\Z_0\to\Z_0}=\norm{C_\varphi}_{\Z\to\Z}=M_\varphi.
\end{equation}
\end{enumerate}
\end{theorem}
\begin{proof}
By Proposition~\ref{prop:conjugacy}, we can work on the sequence spaces. The row functional $x\mapsto(A_\varphi x)(v)$ has norm $\sum_w|a_\varphi(v,w)|$ on both spaces. The triangle inequality gives the upper bound, and for a nonzero row it is attained at the finitely supported sign vector $x_v(w)=\operatorname{sgn}(a_\varphi(v,w))$, with $\operatorname{sgn}(0)=0$. Hence $M_\varphi$ is the exact operator norm whenever $A_\varphi$ is bounded on either space, and $M_\varphi<\infty$ suffices for boundedness on $\ell^\infty(T)$.

Under this row bound, $A_\varphi$ preserves $c_0(T)$ exactly when every column $A_\varphi e_w$ belongs to $c_0(T)$, by density of the finitely supported vectors and closedness of $c_0(T)$ in $\ell^\infty(T)$. An integer-valued column belongs to $c_0(T)$ exactly when its support is finite, since depth balls are finite. This proves~\textup{(b)}.
\end{proof}

\begin{corollary}\label{cor:attainment}
If $M_\varphi<\infty$, then $M_\varphi$ is a positive integer, and $C_\varphi$ attains its $\Z(T)$-operator norm at a function of norm one in $\Z_0(T)$. If $C_\varphi$ preserves the little space, it attains its little-space norm at the same function.
\end{corollary}
\begin{proof}
By~\eqref{eq:G-row}, the root row has sum at least one. The bounded integer row sums therefore attain $M_\varphi$ at some row $v$. The norm-one function $J_0^{-1}x_v$ for the sign vector of that row attains the norm on $\Z(T)$ and, by~\eqref{eq:small-norm}, on $\Z_0(T)$ whenever $C_\varphi$ preserves that space.
\end{proof}

The examples below illustrate the contributions of the root, first-level, and deep rows to the exact norm. Kernel coefficients are combined before taking absolute values.

\begin{example}[Constant symbols]\label{ex:constant}
Let $\varphi(v)=a$ for every $v\in T$. The root row is $G(a,\cdot)$, each first-level row is $-G(a,\cdot)$, and all deeper rows vanish. Thus $C_\varphi$ is bounded on both spaces and
\[
\norm{C_\varphi}_{\Z\to\Z}=\norm{C_\varphi}_{\Z_0\to\Z_0}
=\frac{(|a|+1)(|a|+2)}2.
\]
The fibre $\varphi^{-1}(\{a\})=T$ is infinite, so finite fibres are not necessary for little-space boundedness.
\end{example}

\begin{example}[Cross-branch dilation on a two-ray tree]\label{ex:two-rays}
Let $T=\{o\}\cup\{u_n:n\ge1\}\cup\{v_n:n\ge1\}$, where $u_0=v_0=o$ and the two sequences form rays. For an integer $m\ge1$, define
\[
\varphi_m(o)=o,\qquad\varphi_m(u_n)=v_{mn},\qquad\varphi_m(v_n)=u_{mn}.
\]
The first-order kernel has $m$ coefficients equal to one in each nonroot row, with supports tending to infinity. Thus the norms on $\Lip(T)$ and $\Lip_0(T)$ are both $m$. For the second-order kernel,
\[
(A_{\varphi_m}x)(o)=x(o),\qquad
(A_{\varphi_m}x)(u_1)=(m-1)x(o)+\sum_{j=1}^m(m-j+1)x(v_j).
\]
For $n\ge2$, the three Green rows give a second difference of $P^2x$ over steps of length $m$ along the $v$-ray. Telescoping once for each application of $P$ gives
\begin{equation}\label{eq:two-ray-action}
(A_{\varphi_m}x)(u_n)=\sum_{r=1}^{m}\sum_{s=1}^{m}x(v_{(n-2)m+r+s}).
\end{equation}
The other branch is obtained by interchanging $u,v$. The deeper row sums are $m^2$ and their supports tend to infinity. Since the first-level row sum is $m-1+m(m+1)/2\le m^2$ and the root row sum is $1$,
\[
\norm{C_{\varphi_m}}_{\Z\to\Z}=\norm{C_{\varphi_m}}_{\Z_0\to\Z_0}=m^2.
\]
The maximizing sign vector at $u_2$ in Corollary~\ref{cor:attainment} reconstructs to $\sum_{j=2}^{2m}G_{v_j}$.
\end{example}

\begin{example}[A symbol on the full binary tree]\label{ex:binary}
Let $T=\{0,1\}^{<\infty}$, with the empty word as root. Define $\varphi(w)$ by repeating each letter of $w$ twice and then prepending $0$. Thus $\varphi(o)$ is the one-letter word $0$, and $|\varphi(w)|=2|w|+1$. Formula~\eqref{eq:A-coeff} gives
\[
(A_\varphi x)(o)=2x(o)+x(0),\qquad
(A_\varphi x)(v)=x(\varphi(v)^{--})+2x(\varphi(v)^-)+x(\varphi(v))\quad(v\ne o).
\]
The latter identity includes the first level. Its row sum is $4$, while the root row sum is $3$. The input depths in a nonroot row are $2|v|-1,2|v|,2|v|+1$, so every column has finite support. Hence the norm on both $\Z$ and $\Z_0$ is $4$. The first-order rows have sum $2$ and supports tending to infinity with $|v|$, giving norm $2$ on both Lipschitz spaces.
\end{example}

We now take $T=\Nzero$ with root $0$. The kernel becomes $G(n,k)=(n-k+1)\ind_{\{0\le k\le n\}}$, and the column indexed by $0$ determines the eventual form of a bounded symbol on the little space.

\Needspace{8\baselineskip}
\begin{theorem}\label{thm:affine}
Let $\varphi:\Nzero\to\Nzero$ be arbitrary. Then $C_\varphi$ is bounded on $\Z_0(\Nzero)$ if and only if there are $N,q\in\Nzero$ and $b\in\mathbb Z$ such that
\begin{equation}\label{eq:eventual-affine}
\varphi(n)=qn+b\qquad(n\ge N).
\end{equation}
The eventual slope $q$ is unique. If these conditions hold, then
\begin{equation}\label{eq:affine-norm}
\norm{C_\varphi}_{\Z_0\to\Z_0}
=\norm{C_\varphi}_{\Z\to\Z}
=\max\left\{q^2,\ \max_{0\le n\le N+1}\sum_{k\ge0}|a_\varphi(n,k)|\right\}.
\end{equation}
\end{theorem}
\begin{proof}
By~\eqref{eq:A-coeff} and $G(n,0)=n+1$, the root column satisfies
\[
a_\varphi(n,0)=\varphi(n)-2\varphi(n-1)+\varphi(n-2)\qquad(n\ge2).
\]
Little-space boundedness forces this column to have finite support. Thus the first increments of $\varphi$ eventually equal a fixed integer $q$, giving~\eqref{eq:eventual-affine} with unique slope. Since $\varphi$ takes values in $\Nzero$, $q\ge0$.

Conversely, assume~\eqref{eq:eventual-affine}. If $q=0$, the three images in each row with $n\ge N+2$ coincide, so all these rows vanish. If $q\ge1$, telescoping the two $q$-step differences of $P^2x$ gives, for $n\ge N+2$,
\begin{equation}\label{eq:triangular-row}
\begin{aligned}
(A_\varphi x)(n)
&=\sum_{r=1}^q\bigl[(Px)(\varphi(n-1)+r)-(Px)(\varphi(n-2)+r)\bigr]\\
&=\sum_{r=1}^q\sum_{s=1}^q x(\varphi(n-2)+r+s).
\end{aligned}
\end{equation}
These rows have nonnegative coefficients with sum $q^2$ and support in
$\{\varphi(n-2)+2,\ldots,\varphi(n)\}$. Since $\varphi(n-2)\to\infty$, each fixed column has only finitely many nonzero entries among them.

In both cases, the remaining rows $0\le n\le N+1$ are finite in number and individually finitely supported. They contribute a finite maximum to the row bound and only finitely many entries to each column. Theorem~\ref{thm:boundedness} therefore gives boundedness on both spaces and~\eqref{eq:affine-norm}.
\end{proof}

The next two examples separate boundedness on $\Z$ from preservation of $\Z_0$ and compare the first- and second-order spaces on the chain.

\begin{example}[Large-space boundedness without little-space preservation]\label{ex:alternating}
On $\Nzero$, put $\varphi(n)=0$ for even $n$ and $\varphi(n)=1$ for odd $n$. The kernel rows give
\[
(A_\varphi x)(0)=x(0),\quad (A_\varphi x)(1)=x(1),\quad
(A_\varphi x)(n)=2(-1)^{n+1}\bigl(x(0)+x(1)\bigr)\quad(n\ge2).
\]
Therefore $\norm{C_\varphi}_{\Z\to\Z}=4$. The column indexed by $0$ does not tend to zero, so $C_\varphi$ does not preserve $\Z_0$.
\end{example}

\begin{example}[First-order boundedness does not imply second-order boundedness]\label{ex:floor}
Let $\varphi(n)=\lfloor n/2\rfloor$ on $\Nzero$. The first-order conjugate $DC_\varphi P$ has root row $e_0$, even row $2r$ equal to $e_r$ for $r\ge1$, and zero odd rows. Hence $C_\varphi$ is bounded on both $\Lip$ and $\Lip_0$, with norm $1$.

For the second-order conjugate, the kernel formula gives
\[
a_\varphi(2r,k)=\ind_{\{0\le k\le r\}}\quad(r\ge1),\qquad
a_\varphi(2r+1,k)=-\ind_{\{0\le k\le r\}}\quad(r\ge0),
\]
and $a_\varphi(0,k)=\ind_{\{k=0\}}$. The row sums $r+1$ are unbounded, so Theorem~\ref{thm:boundedness} excludes boundedness on $\Z$. On the little space, the test function $G_0(n)=n+1$ gives $(D^2C_\varphi G_0)(n)=(-1)^n$, so its image lies outside $\Z_0$.
\end{example}

\Needspace{12\baselineskip}
\section{Linear dynamics}\label{sec:dynamics}
A bounded linear operator $A$ on a Banach space $X$ is
\emph{hypercyclic} if some orbit $\{A^n x:n\ge0\}$ is dense in $X$.
On a separable Banach space, hypercyclicity is equivalent to
topological transitivity, which means that for every pair of nonempty
open sets $U,V\subseteq X$, there is an $n\ge0$ with
$A^n(U)\cap V\ne\varnothing$~\cite{LinearChaos}.

A nonzero continuous
linear functional that is bounded on each orbit rules out
hypercyclicity, since the image of a dense orbit would be dense in
$\C$.

We first study Li--Yorke chaos on locally finite rooted trees and then specialize to the chain. Frequent hypercyclicity and Devaney chaos are considered afterward under their respective hypotheses. The notation $\varphi^k$ denotes iteration of the symbol, so $C_\varphi^k=C_{\varphi^k}$. We write
\[
\Per(\varphi)=\{v\in T:\varphi^m(v)=v\text{ for some }m\ge1\},
\qquad \tau_m=\frac{(m+1)(m+2)}2.
\]

\subsection{Li--Yorke chaos}\label{sec:LY}
Two distinct vectors $x,y$ form a \emph{Li--Yorke pair} for $A$ if
\[
\liminf_{n\to\infty}\norm{A^n(x-y)}=0,
\qquad
\limsup_{n\to\infty}\norm{A^n(x-y)}>0.
\]
The operator $A$ is \emph{Li--Yorke chaotic} if it admits an uncountable
\emph{scrambled set}, whose distinct elements form Li--Yorke pairs.
Berm\'udez, Bonilla, Mart\'inez-Gim\'enez, and Peris~\cite{BBMP}
showed that Li--Yorke chaos is equivalent to the existence of an
irregular vector, whose orbit has subsequences with norms tending
to zero and to infinity. In particular, the norms of the powers
must be unbounded.

As in the norm approach to weighted shifts in~\cite[Proposition 3.1 and Theorem 3.2]{KJCK}, we begin with powers. For $C_\varphi\in\mathcal B(\Z_0(T))$, put
\begin{equation}\label{eq:iterate-kernel}
\Gamma_k=\sup_{v\in T}\sum_{w\in T}|a_{\varphi^k}(v,w)|\quad(k\ge0).
\end{equation}
Here $a_{\varphi^k}$ is given by~\eqref{eq:A-coeff}, with all root and first-level conventions unchanged. Theorem~\ref{thm:boundedness} gives
\begin{equation}\label{eq:iterate-norm}
\norm{(\lambda C_\varphi)^k}_{\Z_0\to\Z_0}
=\norm{(\lambda C_\varphi)^k}_{\Z\to\Z}
=|\lambda|^k\Gamma_k.
\end{equation}
Power boundedness therefore excludes Li--Yorke chaos on both spaces.
An irregular vector in $\Z_0(T)$ remains irregular in $\Z(T)$ by the
isometric invariant inclusion. Also, the image sets $\varphi^k(T)$
decrease with $k$. If the finite support of $f$ avoids their
intersection, then $C_\varphi^kf=0$ for all sufficiently large $k$.

We use the dense-stable-set principle of~\cite{BBMP}. If $A^kx\to0$ on a dense subset of a Banach space $X$, then $A$ is Li--Yorke chaotic exactly when $\sup_k\norm{A^k}=\infty$. For sufficiency, the open sets
\[
\bigcup_{k\ge r}\{x:\norm{A^kx}<1/j\},
\qquad
\bigcup_{k\ge0}\{x:\norm{A^kx}>j\}
\quad(j\ge1,\ r\ge0)
\]
are dense. The first family contains the dense stable set. If a set of the second family were not dense, all powers would be bounded by $j$ on an open ball, giving a uniform operator norm bound. By Baire's theorem the intersection is a dense $G_\delta$ set of irregular vectors.

\begin{theorem}\label{thm:LY-tree}
Let $C_\varphi\in\mathcal B(\Z_0(T))$ and assume
\begin{equation}\label{eq:image-escape}
\bigcap_{k\ge0}\varphi^k(T)=\varnothing.
\end{equation}
For every $\lambda\in\C\setminus\{0\}$, the following are equivalent:
\begin{enumerate}[label=\textup{(\roman*)}]
\item $\lambda C_\varphi$ is Li--Yorke chaotic on $\Z_0(T)$;
\item $\lambda C_\varphi$ is Li--Yorke chaotic on $\Z(T)$;
\item $\displaystyle\sup_{k\ge0}|\lambda|^k\Gamma_k=\infty$.
\end{enumerate}
In particular, $C_\varphi$ itself is Li--Yorke chaotic on both spaces.
\end{theorem}
\begin{proof}
Under~\eqref{eq:image-escape}, every finitely supported $f$ satisfies $C_\varphi^kf=0$ eventually. Lemma~\ref{lem:finite-density} gives a dense stable set on $\Z_0(T)$, so the dense-stable-set principle proves (iii)$\Rightarrow$(i). The isometric inclusion and~\eqref{eq:iterate-norm} give (i)$\Rightarrow$(ii)$\Rightarrow$(iii).

For $\lambda=1$, the finite depth balls are eventually avoided, so $|\varphi^k(o)|\to\infty$. By~\eqref{eq:G-row}, the root row gives
\begin{equation}\label{eq:root-iterate-growth}
\Gamma_k\ge\tau_{|\varphi^k(o)|}\longrightarrow\infty.
\end{equation}
This gives the last assertion.
\end{proof}

\begin{example}[Li--Yorke chaos without hypercyclicity on a binary tree]\label{ex:LY-binary}
On $T=\{0,1\}^{<\infty}$, with empty root $o$, put $\varphi(v)=0v$, where $0$ is prepended to the word. Its conjugate acts by
\[
(A_\varphi x)(v)=
\begin{cases}
2x(o)+x(0),&v=o,\\
x(0v)-x(o),&|v|=1,\\
x(0v),&|v|\ge2.
\end{cases}
\]
The rows have uniformly bounded sums and every column has finite support, so $C_\varphi$ is bounded on $\Z_0(T)$. Since $\varphi^k(v)=0^kv$, the absolute row sums of $a_{\varphi^k}$ at the root, first level, and deeper vertices are, respectively,
\[
\tau_k,\qquad \frac{k(k+1)}2+1,\qquad 1.
\]
For the middle formula, the coefficients along $[o,0^k]$ are $-k,-k+1,\ldots,0$, followed by $1$ at $0^kv$. Thus $\Gamma_k=\tau_k$. Moreover, $\inf_v|\varphi^k(v)|=k$, so Theorem~\ref{thm:LY-tree} gives Li--Yorke chaos on both spaces.

Nevertheless, for a fixed $v$ with $|v|\ge2$ and every $k\ge0$,
\[
|(D^2C_\varphi^kf)(v)|=|(D^2f)(0^kv)|\le\norm f_{\Z}.
\]
This nonzero continuous functional is bounded on each orbit, so no orbit can be dense.
\end{example}

We now take $T=\Nzero$ and abbreviate $\Z_0(\Nzero)$ to $\Z_0$. On the chain, periodic vertices may prevent image escape without preventing Li--Yorke chaos. The subspace $Y_F$ from Lemma~\ref{lem:finite-density} isolates this obstruction.

\begin{lemma}\label{lem:periodic-core}
Suppose $C_\varphi\in\mathcal B(\Z_0)$ and $\varphi(n)>n$ for all sufficiently large $n$. Then $F=\Per(\varphi)$ is finite,
\begin{equation}\label{eq:periodic-core}
\bigcap_{k\ge0}\varphi^k(\Nzero)=F,\qquad
\overline{\bigcup_{k\ge1}\ker C_\varphi^k}^{\,\norm{\cdot}_{\Z}}=Y_F.
\end{equation}
Moreover, if $f\in\Z(\Nzero)$ and $\liminf_k\norm{C_\varphi^kf}_{\Z}=0$, then $f|_F=0$.
\end{lemma}
\begin{proof}
Choose $M$ with $\varphi(n)>n$ for $n>M$. Every periodic vertex belongs to all iterated images and lies in $\{0,\ldots,M\}$. Conversely, if $w$ belongs to all iterated images, there are arbitrarily long backward chains ending at $w$. They lie in $\{0,\ldots,\max\{M,w\}\}$ because an orbit starting at a larger vertex increases and cannot reach $w$. A repetition in a sufficiently long chain puts $w$ on a cycle.

Since $\varphi$ permutes $F$, the values of $C_\varphi^kf$ on $F$ are a permutation of the values of $f$ there. If a subsequence tends to zero in norm, continuity of the point evaluations on the finite set $F$ gives $f|_F=0$. In particular, $\ker C_\varphi^k\subseteq Y_F$ for each $k$, and $Y_F$ is closed and invariant. In the other direction, the finite support of $f\in\F_{00}(\Nzero)\cap Y_F$ is disjoint from the intersection of the images. Thus $f$ belongs to some $\ker C_\varphi^k$, and Lemma~\ref{lem:finite-density} gives the asserted closure.
\end{proof}

\begin{theorem}\label{thm:LY-chain}
Let $C_\varphi\in\mathcal B(\Z_0)$, and write $\varphi(n)=qn+b$ for $n\ge N$ as in Theorem~\ref{thm:affine}. The following are equivalent:
\begin{enumerate}[label=\textup{(\roman*)}]
\item $C_\varphi$ is Li--Yorke chaotic on $\Z_0$;
\item $C_\varphi$ is Li--Yorke chaotic on $\Z(\Nzero)$;
\item $\sup_{k\ge0}\Gamma_k=\infty$;
\item $\varphi(n)>n$ for all sufficiently large $n$.
\end{enumerate}
In terms of the affine tail, condition \textup{(iv)} reads $q\ge2$ or $(q=1\text{ and }b>0)$. Neither injectivity nor the absence of periodic vertices is required.
\end{theorem}
\begin{proof}
Assume (iv) and choose $M$ with $\varphi(n)>n$ for $n>M$. Then $F=\Per(\varphi)\subseteq\{0,\ldots,M\}$, and Lemma~\ref{lem:periodic-core} gives a dense stable set in the closed invariant subspace $Y_F$. The fixed test function $G_{M+1}$ vanishes on $F$ and has norm one by~\eqref{eq:column-test}. Its orbit satisfies
\begin{equation}\label{eq:LY-restricted-growth}
\norm{(C_\varphi|_{Y_F})^k}
\ge\frac{(C_\varphi^kG_{M+1})(M+1)}{\tau_{M+1}}
=\frac{\varphi^k(M+1)-M}{\tau_{M+1}}
\longrightarrow\infty,
\end{equation}
by~\eqref{eq:pointbounds}, since $\varphi^k(M+1)$ increases to infinity. The dense-stable-set principle proves (i). The isometric inclusion and~\eqref{eq:iterate-norm} give (i)$\Rightarrow$(ii)$\Rightarrow$(iii).

If (iv) fails, then $q=0$ or $q=1$ with $b\le0$. In the first case, $\varphi(\Nzero)$ is finite and contains all later images. Formula~\eqref{eq:A-coeff} gives $\Gamma_k\le4\tau_{\max\varphi(\Nzero)}$ for $k\ge1$.

In the second case, choose $M\ge N$ with $M\ge\varphi(j)$ for every $j<N$. Since $\varphi(j)=j+b\le j$ for $j\ge N$, the set $\{0,\ldots,M\}$ is forward invariant. If $j>M-kb$, the first $k$ iterates stay above $M$ and equal $j+ib$ at time $i$. If $j\le M-kb$, the orbit starts in this set or reaches it by time $k$. Thus, for every $k\ge0$,
\[
\varphi^k(j)=j+kb\quad(j>M-kb),\qquad
\varphi^k(j)\le M\quad(j\le M-kb).
\]
For a row $n\ge2$ with $n>M-kb+2$, the three images are consecutive, so~\eqref{eq:A-coeff} gives the point mass at $n+kb$, with row sum $1$. For $n\le M-kb+2$, every argument $j$ of $\varphi^k$ used in that row satisfies $j\le M-kb+2$. The displayed formulas imply either $\varphi^k(j)\le M$ or $\varphi^k(j)=j+kb\le M+2$. This also covers the root and first-level rows. The coefficients $1,-2,1$ in~\eqref{eq:A-coeff} and the Green row sum~\eqref{eq:G-row} therefore give
\[
\Gamma_k\le4\tau_{M+2}\qquad(k\ge0).
\]
The integer $M$ was chosen independently of $k$, so this is a uniform bound. Both cases exclude (iii), proving (iii)$\Rightarrow$(iv).
\end{proof}

\begin{example}[Li--Yorke chaos with a fixed root]\label{ex:LY-fixed-root}
Let $\psi(n)=2n$. It induces bounded composition operators on both Zygmund spaces and satisfies $\psi(n)>n$ eventually. Theorem~\ref{thm:LY-chain} gives Li--Yorke chaos, although the fixed root excludes hypercyclicity. Since $\Per(\psi)=\{0\}$, Lemma~\ref{lem:periodic-core} gives $f(0)=g(0)$ for every Li--Yorke pair $(f,g)$. Thus every scrambled set lies in a proper closed affine hyperplane determined by the common value at the root.
\end{example}

The next family quantifies the difference between first- and second-order reconstruction.

\begin{theorem}\label{thm:LY-threshold}
Let $q\ge2$ be an integer and $\varphi_q(n)=q(n+1)-1$. For every $\lambda\ne0$,
\begin{align}
\lambda C_{\varphi_q}\text{ is Li--Yorke chaotic on }\Z_0
&\iff |\lambda|q^2>1,\label{eq:LY-Z-threshold}\\
\lambda C_{\varphi_q}\text{ is Li--Yorke chaotic on }\Lip_0(\Nzero)
&\iff |\lambda|q>1.\label{eq:LY-L-threshold}
\end{align}
The equivalence~\eqref{eq:LY-Z-threshold} also holds on $\Z(\Nzero)$.
\end{theorem}
\begin{proof}
The iterated symbol $\varphi_q^k(n)=q^k(n+1)-1$ is affine with slope $q^k$. In~\eqref{eq:affine-norm} with $N=0$, the root row has sum $\tau_{q^k-1}\le q^{2k}$. The first-level row is nonnegative, since $G(2q^k-1,w)-2G(q^k-1,w)$ equals $w$ for $w<q^k$ and equals $G(2q^k-1,w)$ otherwise. Its sum is therefore $\tau_{2q^k-1}-2\tau_{q^k-1}=q^{2k}$. The first-order rows each contain $q^k$ ones, with the nonroot supports tending to infinity. Thus~\eqref{eq:affine-norm} and~\eqref{eq:L-norm} give
\begin{equation}\label{eq:LY-exact-powers}
\norm{C_{\varphi_q}^k}_{\Z_0\to\Z_0}
=\norm{C_{\varphi_q}^k}_{\Z\to\Z}=q^{2k},
\qquad
\norm{C_{\varphi_q}^k}_{\Lip_0\to\Lip_0}=q^k.
\end{equation}

Since $\min\varphi_q^k(\Nzero)=q^k-1\to\infty$, Theorem~\ref{thm:LY-tree} gives~\eqref{eq:LY-Z-threshold} on both Zygmund spaces. The same escape gives a dense stable set of finitely supported functions in $\Lip_0(\Nzero)$ by Subsection~\ref{subsec:function-spaces}. The dense-stable-set principle now gives~\eqref{eq:LY-L-threshold}.
\end{proof}

\begin{example}[A scalar multiple separating the two spaces]\label{ex:LY-separation}
For $\varphi(n)=2n+1$, formula~\eqref{eq:LY-exact-powers} gives
\[
\norm{(\tfrac13 C_\varphi)^k}_{\Z_0\to\Z_0}=(4/3)^k,\qquad
\norm{(\tfrac13 C_\varphi)^k}_{\Lip_0\to\Lip_0}=(2/3)^k.
\]
Thus $\tfrac13 C_\varphi$ is Li--Yorke chaotic on both Zygmund spaces, while its powers converge to zero in operator norm on $\Lip_0(\Nzero)$.
\end{example}

\subsection{Frequent hypercyclicity}\label{sec:frequent}

A bounded linear operator $A$ is \emph{frequently hypercyclic} if
there is a vector $x$ for which $\{n\ge0:A^n x\in U\}$ has positive
lower density for every nonempty open set $U$~\cite{BayartGrivaux,LinearChaos}.

For $E\subseteq\Nzero$, the lower density is
\[
\underline{\operatorname{dens}}(E)=\liminf_{m\to\infty}
\frac{\#(E\cap\{0,\ldots,m\})}{m+1}.
\]
For composition operators on $L^p$-spaces, chaos and frequent hypercyclicity were studied in~\cite{DP}. We use the following standard form of the Frequent Hypercyclicity Criterion~\cite[Theorem 2.1]{BG}; see also~\cite{BGerr,LinearChaos,LM}.

\begin{criterion}\label{crit:FHC}
Let $A\in\mathcal B(X)$, where $X$ is a separable Banach space. Suppose that there are a dense subset $X_0\subseteq X$ and a map $R:X_0\to X_0$ such that, for each $x\in X_0$,
\[
ARx=x,\qquad \sum_{n=0}^{\infty}A^nx\text{ converges unconditionally},\qquad
\sum_{n=0}^{\infty}R^nx\text{ converges unconditionally}.
\]
Then $A$ is frequently hypercyclic.
\end{criterion}

\Needspace{17\baselineskip}
\begin{theorem}\label{thm:frequent}
Let $\varphi:\Nzero\to\Nzero$ be strictly increasing, with $\varphi(0)>0$, and assume $C_\varphi\in\mathcal B(\Z_0)$. Write $\varphi(n)=qn+b$ for $n\ge N$. The following are equivalent:
\begin{enumerate}[label=\textup{(\roman*)}]
\item $C_\varphi$ satisfies Criterion~\ref{crit:FHC};
\item for every $j\ge1$, there exists $n_j\ge1$ such that
\begin{equation}\label{eq:separation}
\varphi^{n_j}(j)-\varphi^{n_j}(j-1)\ge2;
\end{equation}
\item $q\ge2$;
\item $C_\varphi$ is frequently hypercyclic on $\Z_0$.
\end{enumerate}
If these conditions hold, $n_j=N+1$ works simultaneously for all $j\ge1$.
\end{theorem}

Let $\F_{00}=c_{00}(\Nzero)$ denote the finitely supported functions on the chain, considered as a subspace of $\Z_0$.

\Needspace{18\baselineskip}
\begin{lemma}\label{lem:affine-aux}
Let $\varphi:\Nzero\to\Nzero$ be injective and have no periodic vertices. Suppose $\varphi(n)=qn+b$ for $n\ge N$, where $N,q\in\Nzero$ and $b\in\mathbb Z$. Then $q\ge1$, and the following hold.
\begin{enumerate}[label=\textup{(\roman*)}]
\item If $q=1$, then $b\ge1$. For every $m\ge N+2$, the map $f\mapsto(D^2f)(m)$ is a nonzero continuous functional on $\Z_0$. For $f\in\Z_0$ and $n\ge0$,
\begin{equation}\label{eq:translation-obstruction}
(D^2C_\varphi^nf)(m)=(D^2f)(m+nb)\xrightarrow[n\to\infty]{}0,\qquad
|(D^2C_\varphi^nf)(m)|\le\norm f_{\Z}.
\end{equation}
In particular, $C_\varphi$ has no dense orbit on $\Z_0$.
\item If $q\ge2$, there exists a linear map $R:\F_{00}\to\F_{00}$ such that, for every $f\in\F_{00}$,
\begin{equation}\label{eq:R-properties}
C_\varphi Rf=f,\qquad C_\varphi^nf=0\text{ for all sufficiently large }n,
\qquad\sum_{n=0}^{\infty}\norm{R^nf}_{\Z}<\infty.
\end{equation}
\end{enumerate}
\end{lemma}
\begin{proof}
Theorem~\ref{thm:affine} gives boundedness on $\Z_0$, and injectivity gives $q\ge1$.

For (i), the $N$ distinct initial images must avoid $\{N+b,N+b+1,\ldots\}$, so $N\le N+b$. As $b=0$ would fix the tail, $b\ge1$. For $m\ge N+2$, $\varphi^n$ translates $m,m-1,m-2$ by $nb$. This gives~\eqref{eq:translation-obstruction} by~\eqref{eq:D2} and $D^2f\in c_0$. The functional $f\mapsto(D^2f)(m)$ is continuous, takes value $1$ at $G_m$, and is bounded along each orbit by~\eqref{eq:translation-obstruction}. Thus there is no dense orbit.

For (ii), eventual forward displacement and the absence of periodic vertices give $\bigcap_{n\ge0}\varphi^n(\Nzero)=\varnothing$ by Lemma~\ref{lem:periodic-core}. Hence $C_\varphi^nf=0$ eventually for $f\in\F_{00}$. Since the tail is increasing and $\varphi(n)-n\to\infty$, increase $N\ge2$ so that
\begin{equation}\label{eq:R-tail-choice}
\begin{gathered}
\varphi(n)=qn+b\quad\text{and}\quad \varphi(n)>n\quad(n\ge N-1),\\
\max_{0\le j\le N-2}\varphi(j)<\varphi(N-1).
\end{gathered}
\end{equation}
For $0\le k\le\varphi(N-1)$, set
\begin{equation}\label{eq:R-prefix}
(Rf)(k)=
\begin{cases}
f(j),&k=\varphi(j),\quad0\le j\le N-1,\\
0,&k\notin\varphi(\{0,\ldots,N-1\}).
\end{cases}
\end{equation}
For $k>\varphi(N-1)$, use the unique $j\ge N$ with $\varphi(j-1)<k\le\varphi(j)$ and set
\begin{equation}\label{eq:R-interpolation}
(Rf)(k)=f(j-1)+\frac{k-\varphi(j-1)}q\bigl(f(j)-f(j-1)\bigr).
\end{equation}
Injectivity makes the first rule unambiguous. These linear formulas give $C_\varphi Rf=f$ at the image vertices, and the interpolation values vanish for large $j$ when $f\in\F_{00}$, so $R(\F_{00})\subseteq\F_{00}$.

Put $E_N=\{g\in\F_{00}:g(k)=0\text{ for }k<N\}$. By~\eqref{eq:R-tail-choice}, $\varphi(j)\ge\varphi(N-1)\ge N$ for every $j\ge N-1$. Thus the prefix rule at $k<N$ gives zero or copies the value at the unique predecessor $j<N$ with $\varphi(j)=k$.

If $(R^Nf)(k_0)\ne0$ for some $k_0<N$, repeated use of this rule gives $k_1,\ldots,k_N\in\{0,\ldots,N-1\}$ with $\varphi(k_i)=k_{i-1}$ for $1\le i\le N$. Among these $N+1$ vertices, $k_i=k_j$ for some $0\le i<j\le N$. Then $\varphi^{j-i}(k_j)=k_i=k_j$, contradicting the absence of periodic vertices. Therefore
\begin{equation}\label{eq:R-evacuation}
R^N(\F_{00})\subseteq E_N.
\end{equation}
For $g\in E_N$, $Rg$ vanishes through $\varphi(N-1)$, so $R(E_N)\subseteq E_N$. The first difference on each interpolation interval is
\begin{equation}\label{eq:R-derivative}
(DRg)(k)=\frac{Dg(j)}q
\quad\bigl(\varphi(j-1)<k\le\varphi(j),\ j\ge N\bigr).
\end{equation}
Thus $\norm{Rg}_{\Lip}\le q^{-1}\norm g_{\Lip}$ on $E_N$. Applying this estimate after $N$ iterates gives
\begin{equation}\label{eq:R-geometric}
\norm{R^{N+k}f}_{\Z}\le2\norm{R^{N+k}f}_{\Lip}
\le2q^{-k}\norm{R^Nf}_{\Lip}\quad(k\ge0).
\end{equation}
Summing the geometric bound yields
\begin{equation}\label{eq:R-sum}
\sum_{n=0}^{\infty}\norm{R^nf}_{\Z}
\le\sum_{n=0}^{N-1}\norm{R^nf}_{\Z}
+\frac{2q}{q-1}\norm{R^Nf}_{\Lip}<\infty.
\end{equation}
\end{proof}

\begin{proof}[Proof of Theorem~\ref{thm:frequent}]
Strict increase and $\varphi(0)>0$ give
\begin{equation}\label{eq:forward-growth}
\varphi(j)\ge j+1,\qquad\varphi^n(j)\ge j+n.
\end{equation}
Thus $\varphi$ is injective and has no periodic vertices. For $q\ge2$, Lemma~\ref{lem:affine-aux}\textup{(ii)} verifies Criterion~\ref{crit:FHC} on the dense subspace $\F_{00}$, giving (iii)$\Rightarrow$(i)$\Rightarrow$(iv). By Lemma~\ref{lem:affine-aux}\textup{(i)}, (iv) excludes $q=1$, so (iv)$\Rightarrow$(iii).

For $j\ge1$, both $\varphi^N(j-1)$ and $\varphi^N(j)$ lie in the affine tail by~\eqref{eq:forward-growth}. They are distinct and ordered, so if $q\ge2$,
\[
\varphi^{N+1}(j)-\varphi^{N+1}(j-1)
=q\bigl(\varphi^N(j)-\varphi^N(j-1)\bigr)\ge2.
\]
If $q=1$, then $b\ge1$ and $\varphi^n(N+1)-\varphi^n(N)=1$ for every $n$. This proves (ii)$\Leftrightarrow$(iii) and the asserted choice of $n_j$.
\end{proof}

The separation condition parallels the little Lipschitz characterization in~\cite[Theorem 3.1]{LM}. Under $J_0$, the criterion for the sequence conjugate uses $J_0RJ_0^{-1}$ on the dense subspace $J_0(\F_{00})\subseteq c_0$.

\begin{example}[An expanding affine symbol]\label{ex:frequent-affine}
Let $\varphi(n)=2n+1$. The Green rows give
\begin{equation}\label{eq:2n1-matrix}
(A_\varphi x)(0)=2x(0)+x(1),\qquad
(A_\varphi x)(n)=x(2n-1)+2x(2n)+x(2n+1)\quad(n\ge1).
\end{equation}
Every column has finite support. The root row sum is $3$, and every other row sum is $4$. Therefore
\[
\norm{C_\varphi}_{\Z_0\to\Z_0}=4.
\]
The eventual slope is $2$, so Theorem~\ref{thm:frequent} gives the Frequent Hypercyclicity Criterion and frequent hypercyclicity. At $n=1$, the maximizing sign vector in Corollary~\ref{cor:attainment} reconstructs to $G_1+G_2+G_3$.
\end{example}

\begin{example}[A bounded translation without frequent hypercyclicity]\label{ex:translation}
Let $\psi(n)=n+1$. Its conjugate is
\[
(A_\psi x)(0)=2x(0)+x(1),\quad
(A_\psi x)(1)=-x(0)+x(2),\quad
(A_\psi x)(n)=x(n+1)\quad(n\ge2).
\]
The matrix preserves $c_0$ and has norm $3$, so $C_\psi$ is bounded on $\Z_0$ with the same norm. Theorem~\ref{thm:LY-chain} gives Li--Yorke chaos on both Zygmund spaces. Lemma~\ref{lem:affine-aux}\textup{(i)} with $b=1$ rules out hypercyclicity via~\eqref{eq:translation-obstruction}, so the Frequent Hypercyclicity Criterion also fails.
For comparison, the operators in this example and Example~\ref{ex:frequent-affine} both have $\Lip_0$-operator norm $2$, as follows from their first-order kernel rows.
\end{example}

\Needspace{14\baselineskip}
\subsection{Mixing and Devaney chaos}\label{sec:chaos}
A bounded linear operator $A$ on a Banach space $X$ is
\emph{topologically mixing} if $A^n(U)\cap V\ne\varnothing$
for every pair of nonempty open sets $U,V\subseteq X$ and all
sufficiently large $n$. It is \emph{Devaney chaotic} if it is
hypercyclic and its periodic vectors
$\Per(A)=\bigcup_{m\ge1}\ker(A^m-I)$ are dense in $X$~\cite{LinearChaos}.

For related periodic-point methods for shifts on trees, see~\cite{GPchaos}.

\begin{theorem}\label{thm:chaos}
Let $\varphi:\Nzero\to\Nzero$ and suppose $C_\varphi\in\mathcal B(\Z_0)$. Let $q\in\Nzero$ be the eventual slope from Theorem~\ref{thm:affine}. The following are equivalent:
\begin{enumerate}[label=\textup{(\roman*)}]
\item $C_\varphi$ is Devaney chaotic on $\Z_0$;
\item $C_\varphi$ is topologically mixing on $\Z_0$;
\item $\varphi$ is injective, has no periodic vertices, and $q\ge2$.
\end{enumerate}
\end{theorem}
\begin{proof}
Both (i) and (ii) imply hypercyclicity on the separable space $\Z_0$. If $\varphi(u)=\varphi(v)$ for $u\ne v$, the nonzero continuous functional $f\mapsto f(u)-f(v)$ vanishes on each positive iterate. If $v$ is periodic, evaluation at $v$ takes finitely many values along each orbit. Both cases preclude a dense orbit, so hypercyclicity requires injectivity and absence of periodic vertices. Injectivity excludes $q=0$, and Lemma~\ref{lem:affine-aux}\textup{(i)} excludes $q=1$. This proves necessity of (iii).

Assume (iii) and use $R$ from Lemma~\ref{lem:affine-aux}\textup{(ii)}. For $u,v\in\F_{00}$, $u+R^nv\to u$ and, for all sufficiently large $n$,
\[
C_\varphi^n(u+R^nv)=v
\]
by~\eqref{eq:R-properties}. Density of $\F_{00}$ gives mixing, by the right-inverse argument of Kitai's Criterion~\cite[Section 1]{GG}.

For $f\in\F_{00}$ and $m\ge1$, the series
\begin{equation}\label{eq:periodic-series}
p_{m,f}=\sum_{j=0}^{\infty}R^{jm}f
\end{equation}
converges absolutely by~\eqref{eq:R-sum}. For $m$ large enough that $C_\varphi^mf=0$, termwise application of $C_\varphi^m$ and $C_\varphi R=I$ give
\[
C_\varphi^mp_{m,f}=\sum_{j=1}^{\infty}R^{(j-1)m}f=p_{m,f}.
\] Since
\begin{equation}\label{eq:periodic-approximation}
\norm{p_{m,f}-f}_{\Z}\le\sum_{n\ge m}\norm{R^nf}_{\Z}\longrightarrow0,
\end{equation}
periodic vectors are dense. Together with mixing, this gives Devaney chaos.
\end{proof}

\begin{example}[Global affine symbols and explicit periodic vectors]\label{ex:global-chaos}
For $\varphi(n)=qn+b$ with $q,b\in\Nzero$, Theorem~\ref{thm:affine} gives boundedness on $\Z_0$. By Theorem~\ref{thm:chaos}, the operator is Devaney chaotic exactly when $q\ge2$ and $b\ge1$. For $q\ge2$, $b=0$ fixes the root, while $b\ge1$ gives $\varphi(n)>n$ at every vertex.

For $q\ge2$ and $b\ge1$, set $c=b/(q-1)>0$. Given $m\ge1$ and $\ell\in\mathbb Z$, put
\[
\theta=\frac{2\pi\ell}{m\log q},\qquad
\rho_{m,\ell}(n)=\exp\bigl(i\theta\log(n+c)\bigr).
\]
Since $\varphi(n)+c=q(n+c)$, one has $C_\varphi^m\rho_{m,\ell}=\rho_{m,\ell}$. The second derivative of $t\mapsto(t+c)^{i\theta}$ is $i\theta(i\theta-1)(t+c)^{i\theta-2}$, which tends to zero as $t\to\infty$. The second finite differences therefore tend to zero, so $\rho_{m,\ell}\in\Z_0$. Density of periodic vectors follows from~\eqref{eq:periodic-series}--\eqref{eq:periodic-approximation}.
\end{example}

\begin{example}[Devaney chaos for a nonmonotone symbol]\label{ex:nonmonotone}
Define
\[
\varphi(0)=2,\qquad \varphi(1)=0,\qquad\varphi(n)=2n\quad(n\ge2).
\]
The symbol has an affine tail of slope $2$, so $C_\varphi$ is bounded on $\Z_0$ by Theorem~\ref{thm:affine}. It is injective, and its initial orbit is $1\mapsto0\mapsto2\mapsto4\mapsto8\mapsto\cdots$. Every orbit entering $\{2,3,\ldots\}$ moves strictly to the right, so there are no periodic vertices. Theorem~\ref{thm:chaos} gives Devaney chaos and mixing. Lemmas~\ref{lem:finite-density} and~\ref{lem:affine-aux}\textup{(ii)} also verify Criterion~\ref{crit:FHC} directly, giving frequent hypercyclicity for this nonmonotone symbol.

The row sums for $n=0,1,2,3$ are, respectively, $6,11,19,9$, and every row with $n\ge4$ has sum $4$. For instance,
\[
(A_\varphi x)(2)=6x(0)+6x(1)+4x(2)+2x(3)+x(4).
\]
Consequently $\norm{C_\varphi}_{\Z_0\to\Z_0}=19$. The finite initial segment changes the exact norm substantially without preventing the dynamical properties above.
\end{example}

\begin{example}[An expansive tail with a periodic obstruction]\label{ex:cycle}
Let $\psi(0)=1$, $\psi(1)=0$, and $\psi(n)=2n$ for $n\ge2$. The symbol is injective and eventually affine with slope $2$, so its composition operator is bounded on $\Z_0$. However, $0$ and $1$ form a two-cycle: $(C_\psi^nf)(0)\in\{f(0),f(1)\}$ for all $n\ge0$. No orbit can therefore be dense in $\Z_0$, and the operator is neither hypercyclic nor Devaney chaotic. It is nevertheless Li--Yorke chaotic by Theorem~\ref{thm:LY-chain}. Thus the absence of periodic vertices in Theorem~\ref{thm:chaos} cannot be omitted.
\end{example}

We recall the standard \emph{Hypercyclicity Criterion}~\cite[Chapter 3]{LinearChaos}.
An operator $A\in\mathcal B(X)$ on a separable Banach space satisfies
this criterion if there exist dense subsets $X_0,Y_0\subseteq X$,
a strictly increasing sequence $(n_k)$ of positive integers, and maps
$S_k:Y_0\to X$ such that, as $k\to\infty$,
\[
A^{n_k}x\longrightarrow0,\qquad
S_ky\longrightarrow0,\qquad
A^{n_k}S_ky\longrightarrow y
\quad(x\in X_0,\ y\in Y_0).
\]
Every operator satisfying this criterion is hypercyclic.

We conclude by relating hypercyclicity, the Hypercyclicity Criterion,
and Devaney chaos under the strictly increasing hypotheses of
Theorem~\ref{thm:frequent}. By Theorem~\ref{thm:LY-chain}, Li--Yorke
chaos already holds under these hypotheses, including the translation case.

\Needspace{12\baselineskip}
\begin{corollary}\label{cor:increasing-chaos}
Let $\varphi:\Nzero\to\Nzero$ be strictly increasing, assume $\varphi(0)>0$, and suppose $C_\varphi\in\mathcal B(\Z_0)$. Write $\varphi(n)=qn+b$ for $n\ge N$, with $q\in\N$. The following are equivalent:
\begin{enumerate}[label=\textup{(\roman*)}]
\item $C_\varphi$ is hypercyclic;
\item $C_\varphi$ satisfies the Hypercyclicity Criterion;
\item $C_\varphi$ is Devaney chaotic;
\item $q\ge2$.
\end{enumerate}
Under these assumptions, the four conditions are also equivalent to frequent hypercyclicity and to satisfaction of the Frequent Hypercyclicity Criterion.
\end{corollary}
\begin{proof}
Strict increase and~\eqref{eq:forward-growth} give injectivity and absence of periodic vertices. Hence (iii)$\Leftrightarrow$(iv) follows from Theorem~\ref{thm:chaos}.

Assume (iv), and let $R$ be the right inverse in Lemma~\ref{lem:affine-aux}\textup{(ii)}. For every $f\in\F_{00}$, \eqref{eq:R-properties} gives
\[
C_\varphi^nf\longrightarrow0,\qquad
R^nf\longrightarrow0,\qquad
C_\varphi^nR^nf=f.
\]
The second limit follows from the summability of $\norm{R^nf}_{\Z}$, and the last identity follows by iterating $C_\varphi R=I$ on $\F_{00}$. Taking $X_0=Y_0=\F_{00}$, $n_k=k$, and $S_k=R^k$ verifies the Hypercyclicity Criterion. This proves (iv)$\Rightarrow$(ii)$\Rightarrow$(i). The necessity argument in the proof of Theorem~\ref{thm:chaos} gives (i)$\Rightarrow$(iv). The final assertion follows from Theorem~\ref{thm:frequent}.
\end{proof}

\end{document}